\documentstyle[12pt,amssymb]{article}
\begin{document}

{\Large

\noindent{\bf A set of orthogonal polynomials, dual to}

\noindent{\bf alternative $q$-Charlier polynomials} }
\bigskip

\noindent{\sc Natig M. Atakishiyev}

\noindent {\it Instituto de Matem\'aticas, UNAM, CP 62210
Cuernavaca, Morelos, M\'exico}
 \medskip

\noindent and
 \medskip

\noindent{\sc Anatoliy U. Klimyk}

\noindent {\it Bogolyubov Institute for Theoretical Physics, 03143
Kiev, Ukraine}

\medskip

 \noindent
E-mail: natig@matcuer.unam.mx and aklimyk@bitp.kiev.ua
\medskip

\begin{abstract}
The aim of this paper is to derive (by using two operators,
representable by a Jacobi matrix) a family of $q$-orthogonal
polynomials, which turn to be dual to alternative $q$-Charlier
polynomials. A discrete orthogonality relation and a three-term
recurrence relation for these dual polynomials are explicitly
obtained. The completeness property of dual alternative
$q$-Charlier polynomials is also established.
\end{abstract}

\bigskip


\medskip
\bigskip

\noindent{\bf 1. Introduction}
\bigskip

It is well known that orthogonal polynomials are closely connected
with spectral properties of symmetric operators, which can be
represented in some basis by a Jacobi matrix, and with the
classical moment problem. Namely, a spectrum of an operator,
represented by a Jacobi matrix, is determined by an orthogonality
measure for corresponding orthogonal polynomials. If orthogonal
polynomials admit many orthogonality relations, then the
corresponding symmetric operator is not self-adjoint and it leaves
room for infinitely many self-adjoint extensions. These extensions
are determined by orthogonality measures for the appropriate
orthogonal polynomials.

Contrary to orthogonal polynomials of the hypergeometric type
(Wilson, Jacobi, Laguerre and so on), basic hypergeometric
polynomials (or $q$-orthogonal polynomials) are not so deeply
understood yet. These polynomials are collected in the
$q$-analogue of the Askey--scheme of orthogonal polynomials (see,
for example, [1]), which starts with Askey--Wilson polynomials and
$q$-Racah polynomials, introduced in [2] and [3]. Importance of
these polynomials is magnified by the fact that they are closely
related to the theory of quantum groups. As an instance of such
connection we refer to a paper [4], in which Al-Salam--Chihara
$q$-orthogonal polynomials have been employed to construct locally
compact quantum group $SU_q(1,1)$. Another application of
$q$-orthogonal polynomials is related to the theory of
$q$-difference equations, which often surface in contemporary
theoretical and mathematical physics.

The purpose of present paper is to study completeness and duality
properties of alternative $q$-Charlier polynomials. We shall show
below that this originates a novel type of $q$-orthogonal
polynomials and a discrete orthogonality relation for them. To
achieve this, we essentially use two operators $I_1$ and
$q^{J_0}$, which are certain representation operators for the
quantum algebra $U_q({\rm su}_{1,1})$ with a lowest weight
(however, we do not use explicitly the theory of representations
in what follows). The operator $I_1$ is related to the three-term
recurrence relation for alternative $q$-Charlier polynomials. We
diagonalize the trace class operator $I_1$ and obtain two bases in
the Hilbert space: an initial basis and a basis of normalized
eigenvectors of $I_1$. These bases are connected by an orthogonal
matrix. The orthogonality relations for rows and columns of this
matrix lead to orthogonality relations for alternative
$q$-Charlier polynomials and for the functions, which are dual to
these polynomials. We extract from the latter functions a dual set
of polynomials and obtain a discrete orthogonality relation for
them. As a result, one is led to the completeness property of dual
alternative $q$-Charlier polynomials.

Observe that the present paper is a continuation of our research,
initiated in [5] and [6].

Throughout the sequel we always assume that $q$ is a fixed
positive number such that $q<1$. We use (without additional
explanation) notations of the theory of special functions and the
standard $q$-analysis (see, for example, [7] or [8]).
\bigskip

\noindent{\bf 2. Pair of operators $(I_1,J)$}
\bigskip

Let ${\cal H}$ be a separable complex Hilbert space with an
orthonormal basis $|n\rangle$, $n = 0,1,2,\cdots$. We define on
${\cal H}$ two operators. The first one, denoted as $q^{J_0}$,
acts on the basis elements as $q^{J_0}|n\rangle =q^n |n\rangle$,
and the second one, denoted as $I_1$, is given by the formula
$$
I_1 |n\rangle =a_n|n+1\rangle +a_{n-1}|n-1\rangle +b_n|n\rangle  ,
\eqno (1)
$$
$$
a_n=-(aq^{3n+1})^{1/2}\frac{ \sqrt{(1-q^{n+1})(1+aq^{n})}}
{(1+aq^{2n+1})\sqrt{(1+aq^{2n}) (1+aq^{2n+2})}} ,
$$  $$
b_n=q^n\left( \frac{1+aq^n}{(1+aq^{2n}) (1+aq^{2n+1})}+aq^{n-1}
\frac{1-q^{n}}{(1+aq^{2n-1}) (1+aq^{2n})}  \right) ,
$$
where $a$ is any fixed positive number. Clearly, $I_1$ is a
symmetric operator.

Since $a_n\to 0$ and $b_n\to 0$ when $n\to \infty$, the operator
$I_1$ is bounded. Therefore, we assume that it is defined on the
whole Hilbert space ${\cal H}$. For this reason, $I_1$ is a
self-adjoint operator. Let us show that $I_1$ is a trace class
operator. For the coefficients $a_n$ and $b_n$ from (1), we have
$a_{n+1}/a_n \to q^{3/2}$ and $b_{n+1}/b_n \to q$ when $n\to
\infty$.  Since $0<q<1$, for the sum of all matrix elements of the
operator $I_1$ in the basis $|n\rangle$, $n = 0,1,2,\cdots$, we
have $\sum _n (2a_n+b_n)< \infty$. This means that $I_1$ is a
trace class operator. Thus, a spectrum of $I_1$ is discrete and
has a single accumulation point at 0. Moreover, a spectrum of
$I_1$ is simple, since $I_1$ is representable by a Jacobi matrix
with $a_n\ne 0$ (see [9], Chapter VII).

To find eigenfunctions $\xi_\lambda $ of the operator $I_1$, $I_1
\xi_\lambda =\lambda \xi_\lambda $, we set $\xi_\lambda =\sum _n
\beta_n(\lambda)|n \rangle$, where $\beta_n(\lambda)$ are
appropriate numerical coefficients. Acting by the operator $I_1$
upon both sides of this relation, one derives that $\sum
_{n=0}^{\infty}\, \beta_n(\lambda)\, (a_n|n+1\rangle
+a_{n-1}|n-1\rangle +b_n|n\rangle )= \lambda \sum_{n=0}^{\infty}\,
\beta_n(\lambda) |n\rangle$, where $a_n$ and $b_n$ are the same as
in (1). Collecting in this identity all factors, which multiply
$|n\rangle$ with fixed $n$, one derives the recurrence relation
for the coefficients $\beta_n(\lambda)$:
$$
\beta_{n+1}(\lambda)a_n +\beta_{n-1}(\lambda)a_{n-1}+
\beta_{n}(\lambda)b_n= \lambda \beta_{n}(\lambda).
$$
The substitution
$$
\beta_{n}(\lambda)=\left( \frac{(-a;q)_n\,(1+aq^{2n})} {(q;q)_n\,
(1+a)(a/q)^n}\right) ^{1/2} q^{-n(n+3)/4} \beta'_{n}(\lambda)
$$
reduces this relation to the following one
$$
-A_n \beta'_{n+1}(\lambda)- C_n \beta'_{n-1}(\lambda)
+(A_n+C_n)\beta'_{n}(\lambda)=\lambda \beta'_{n}(\lambda),
$$
$$
A_n=q^n\frac{1+aq^{n}}{(1+aq^{2n}) (1+aq^{2n+1})}, \ \ \ \
C_n=aq^{2n-1}\frac{1-q^{n}}{(1+aq^{2n-1}) (1+aq^{2n})}.
$$
This is the recurrence relation for the alternative $q$-Charlier
polynomials
$$
K_n(\lambda ;a;q):={}_2\phi_1 (q^{-n}, -aq^{n};\; 0; \; q,q\lambda
)
$$
(see, formulas (3.22.1) and (3.22.2) in [1]). Therefore,
$\beta'_n(\lambda )=K_n(\lambda ;a;q)$ and
$$
\beta_n(\lambda )=\left( \frac{(-a;q)_n\,(1+aq^{2n})} {(q;q)_n\,
(1+a)a^n}\right) ^{1/2} q^{-n(n+1)/4} K_n(\lambda ;a;q). \eqno (2)
$$
For the eigenvectors $\xi _\lambda$ we thus have the expression
$$
\xi _\lambda =\sum_{n=0}^\infty\, \left(
\frac{(-a;q)_n\,(1+aq^{2n})} {(q;q)_n\, (1+a)a^n}\right) ^{1/2}
q^{-n(n+1)/4} K_n(\lambda ;a;q) |n\rangle .
 \eqno (3)
$$
Since the spectrum of the operator $I_1$ is discrete, only for a
discrete set of values of $\lambda$ these vectors belong to the
Hilbert space ${\cal H}$. This discrete set of eigenvectors
determines a spectrum of $I_1$.

Now we look for a spectrum of the operator $I_1$ and for a set of
polynomials, dual to alternative $q$-Charlier polynomials. To this
end we use the action of the operator
$$
J:= q^{-J_0} - a\,q^{J_0}
$$
upon the eigenvectors $\xi _\lambda$, which belong to the Hilbert
space ${\cal H}$. In order to find how this operator acts upon
these vectors, one can use the $q$-difference equation
$$
(q^{-n}-aq^{n})K_n(\lambda)=- aK_{n} (q\lambda
)+\lambda^{-1}K_n(\lambda)-\lambda^{-1}
(1-\lambda)K_n(q^{-1}\lambda) \eqno(4)
$$
for the alternative $q$-Charlier polynomials $K_n(\lambda)\equiv
K_n(\lambda ;a;q)$ (see formula (3.22.5) in [1]). Multiply both
sides of (4) by $d_n\,|n\rangle$ and sum up over $n$, where $d_n$
are the coefficients of $K_n(\lambda ;a;q)$ in the expression (2)
for $\beta_n(\lambda)$. Taking into account the formula (3) and
the fact that $J|n\rangle=(q^{-n}-aq^{n})|n\rangle$, one obtains
the relation
$$
J\,\xi _{\lambda}= -a\,\xi _{q\lambda}+ \lambda^{-1}\, \xi
_{\lambda}- \lambda^{-1}(1-\lambda)\, \xi_{q^{-1}\lambda}. \eqno
(5)
$$

We shall see in the next section that the spectrum of the operator
$I_1$ consists of the points $q^n$, $n=0,1,2,\cdots$. This means
that $J$ has the form of a Jacobi matrix in the basis of
eigenvectors of $I_1$; that is, the pair of the operators $I_1$
and $J$ form a Leonard pair (see [9] for the corresponding
definition).
\bigskip

\noindent{\bf 3. Spectrum of $I_1$ and orthogonality of
alternative $q$-Charlier polynomials}
\medskip

The aim of this section is to find, by using the operators $I_1$
and $J$, a basis in the Hilbert space ${\cal H}$, which consists
of eigenvectors of the operator $I_1$ in a normalized form, and to
derive explicitly the unitary matrix $U$, connecting this basis
with the basis $|n\rangle$, $n=0,1,2,\cdots$, in ${\cal H}$. This
matrix leads directly to the orthogonality relation for
alternative $q$-Charlier polynomials. For this purpose we first
find a spectrum of $I_1$.

Let us analyze a form of the spectrum of $I_1$ from the point of
view of the spectral theory of trace class operators. If $\lambda$
is a spectral point of the operator $I_1$, then (as it is easy to
see from (5)) a successive action by the operator $J$ upon the
vector (eigenvector of $I_1$) $\xi_\lambda$ leads to the
eigenvectors $\xi_{q^m\lambda}$, $m=0,\pm 1, \pm 2,\cdots$.
However, since $I_1$ is a trace class operator, not all of these
points may belong to the spectrum of $I_1$, since $q^{-m}\lambda
\to\infty$ when $m\to +\infty$ if $\lambda\ne 0$. This means that
the coefficient $1-\lambda' $ of $\xi _{q^{-1}\lambda'}$ in (5)
must vanish for some eigenvalue $\lambda'$. Clearly, it vanishes
when $\lambda' =1$. Moreover, this is the only possibility for the
coefficient of $\xi _{q^{-1}\lambda'}$ in (5) to vanish, that is,
the point $\lambda =1$ is a spectral point for the operator $I_1$.
Let us show that the corresponding eigenfunction $\xi _{1}\equiv
\xi_{q^{0}}$ belongs to the Hilbert space ${\cal H}$.

By formula (II.6) of Appendix II in [7], one has $K_n(1
;a;q)={}_2\phi_1 (q^{-n}, -aq^{n};\; 0; \; q,q )=(-a)^nq^{n^2}$.
Therefore,
$$
\langle \xi_1,\xi_1\rangle = \sum_{n=0}^\infty
\frac{(-a;q)_n(1+aq^{2n})}{(1+a)(q;q)_na^nq^{n(n+1)/2}}K^2_n(1
;a;q) = \sum_{n=0}^\infty
\frac{(-a;q)_n(1+aq^{2n})a^n}{(1+q)(q;q)_n q^{-n(3n-1)/2}} .
 \eqno (6)
$$
In order to calculate this sum, we take the limit $d,e\to \infty$
in the equality
$$
 \sum_{n=0}^\infty
\frac{(1+aq^{2n})(-a;q)_n(d;q)_n(e;q)_n}{(1+a)(-aq/d;q)_n(-aq/e;q)_n(q;q)_n}
 \left( \frac{aq}{de}\right) ^n q^{n(n-1)/2} =\frac{(-aq;q)_\infty
 (-aq/de;q)_\infty}{(-aq/d;q)_\infty (-aq/e;q)_\infty}
 $$
(see formula in Exercise 2.12, Chapter 2 of [7]).  Since
$$
\lim_{d,e\to \infty}(d;q)_n(e;q)_n (aq/de)^n=q^{n(n-1)}(aq)^n ,
$$
we obtain from here that the sum in (6) is equal to $
(-aq;q)_\infty$, that is, $\langle \xi_1,\xi_1\rangle <\infty$ and
$\xi_1$ belongs to the Hilbert space ${\cal H}$. Thus, the point
$\lambda =1$ does belong to the spectrum of the operator $I_1$.

Let us find other spectral points of the operator $I_1$ (recall
that a spectrum of $I_1$ is discrete). Setting $\lambda = 1$ in
(5), we see that the operator $J$ transforms $\xi _{q^0}$ into a
linear combination of the vectors $\xi _{q}$ and $\xi_{q^0}$.
Moreover, $\xi_q$ belongs to the Hilbert space ${\cal H}$, since
the series
$$
\langle \xi_q,\xi_q\rangle = \sum_{n=0}^\infty
\frac{(-a;q)_n\,(1+aq^{2n})}{(1+a)(q;q)_n\,a^n}q^{-n(n+1)/2}\,K^2_n(q
;a;q)
$$
is majorized by the corresponding series (6) for $\xi_{q^0}$.
Therefore, $\xi _{q}$ belongs to the Hilbert space ${\cal H}$ and
the point $q$ is an eigenvalue of the operator $I_1$. Similarly,
setting $\lambda=q$ in (5), we find that $\xi _{q^2}$ is an
eigenvector of $I_1$ and the point $q^2$ belongs to the spectrum
of $I_1$. Repeating this procedure, we find that all $\xi _{q^n}$,
$n=0,1,2,\cdots$, are eigenvectors of $I_1$ and the set $q^n$,
$n=0,1,2,\cdots$, belongs to the spectrum of $I_1$. So far, we do
not know yet whether other spectral points exist or not.

The vectors $\xi _{q^n}$, $n=0,1,2,\cdots$, are linearly
independent elements of the Hilbert space ${\cal H}$ (since they
correspond to different eigenvalues of the self-adjoint operator
$I_1$). Suppose that values $q^n$, $n=0,1,2,\cdots$, constitute a
whole spectrum of $I_1$. Then the set of vectors $\xi _{q^n}$,
$n=0,1,2,\cdots$, is a basis in the Hilbert space ${\cal H}$.
Introducing the notation $\Xi _k:=\xi_{q^k}$, $k=0,1,2,\cdots$, we
find from (5) that
$$
J \,\Xi _k = - a \,\Xi _{k+1} + q^{-k}\, \Xi _k - q^{-k}(1-q^k)\,
\Xi _{k-1} .
$$
As we see, the matrix of the operator $J$ in the basis $\Xi _k$,
$k=0,1,2,\cdots$, is not symmetric, although in the initial basis
$|n\rangle$, $n=0,1,2,\cdots$, it was symmetric. The reason is
that the matrix $(a_{mn})$ with entries $a_{mn}:=\beta_m(q^n)$,
$m,n=0,1,2,\cdots$, where $\beta_m(q^n)$ are the coefficients (2)
in the expansion $\xi _{q^n}=\sum _m \,\beta_m(q^n)|n\rangle$, is
not unitary. This fact is equivalent to the statement that the
basis $\Xi _n=\xi_{q^n}$, $n=0,1,2,\cdots$, is not normalized. To
normalize it, one has to multiply $\Xi _n$ by corresponding
numbers $c_n$ (which are not known at this moment). Let $\hat\Xi
_n = c_n\Xi _n$, $n=0,1,2,\cdots$, be a normalized basis. Then the
matrix of the operator $J$ is symmetric in this basis. Since  $J$
has in the basis $\{ \hat\Xi _n\}$ the form
$$
J\, \hat\Xi _n = -c_{n+1}^{-1}c_n a\, \hat\Xi _{n+1} + q^{-n}\,
\hat\Xi _n - c_{n-1}^{-1}c_n q^{-n}(1-q^n) \,\hat\Xi_{n-1} ,
$$
then its symmetricity means that
$c_{n+1}^{-1}c_na=c_{n}^{-1}c_{n+1} q^{-n-1}(1-q^{n+1})$, that is,
$c_{n}/c_{n-1} =\sqrt{aq^n/(1-q^n)}$. Therefore,
$$
c_n= c(a^nq^{n(n+1)/2}/(q;q)_n )^{1/2},
$$
where $c$ is a constant.

The expansions
$$
\hat\xi _{q^n}(x)\equiv \hat\Xi _n(x)= \sum _m
c_n\beta_m(q^n)|m\rangle \equiv \sum _m {\hat a}_{mn}|m\rangle
\eqno (7)
$$
connect two orthonormal bases in the Hilbert space ${\cal H}_l$.
This means that the matrix $({\hat a}_{mn})$, $m,n=0,1,2,\cdots$,
with entries
$$
{\hat a}_{mn}=c_n\beta _m(q^n)= c\left(
\frac{a^nq^{n(n+1)/2}}{(q;q)_n}
\frac{(-a;q)_m\,(1+aq^{2m})}{(1+a)(q;q)_m\,a^m q^{m(m+1)/2}}
\right)^{1/2} K_m(q^n ;a;q) \eqno(8)
$$
is unitary, provided  that the constant $c$ is appropriately
chosen. In order to calculate this constant, we use the relation
$\sum_{m=0}^\infty |{\hat a}_{mn}|^2=1$ for $n=0$. Then this sum
is a multiple of the sum in (6) and, consequently,
$c=(-aq;q)^{-1/2}_\infty$.

The matrix $({\hat a}_{mn})$ is real and orthogonal, that is,
$$   {\textstyle
\sum _n {\hat a}_{mn}{\hat a}_{m'n}=\delta_{mm'},\ \ \ \ \sum _m
{\hat a}_{mn}{\hat a}_{mn'}=\delta_{nn'} .  } \eqno (9)
$$
Substituting into the first sum over $n$ in (9) the expressions
for ${\hat a}_{mn}$, we obtain the identity
$$
\sum_{n=0}^\infty
\frac{a^nq^{n(n+1)/2}}{(q;q)_n}\,K_m(q^n;a;q)\,K_{m'}(q^n;a;q) =
\frac{(-aq^m;q)_\infty \, a^m\,(q;q)_m} {(1+aq^{2m})}\,
q^{m(m+1)/2} \delta_{mm'}\,, \eqno (10)
$$
which must yield the orthogonality relation for alternative
$q$-Charlier polynomials. An only gap, which remains to be
clarified, is the following. We have assumed that the points
$q^n$, $n=0,1,2,\cdots$, exhaust the whole spectrum of $I_1$. Let
us show that this is the case.

Recall that the self-adjoint operator $I_1$ is represented by a
Jacobi matrix in the basis $|n\rangle$, $n=0,1,2,\cdots$.
According to the theory of operators of such type (see, for
example, [9], Chapter VII), eigenvectors $\xi_\lambda$ of $I_1$
are expanded into series in the basis $|n\rangle$,
$n=0,1,2,\cdots$, with coefficients, which are polynomials in
$\lambda$. These polynomials are orthogonal with respect to some
positive measure $d\mu (\lambda)$ (moreover, for self-adjoint
operators this measure is unique). The set (a subset of ${\Bbb
R}$), on which these polynomials are orthogonal, coincides with
the spectrum of the operator under consideration and the spectrum
is simple.

We have found that the spectrum of $I_1$ contains the points
$q^n$, $n=0,1,2,\cdots$. If the operator $I_1$ would have other
spectral points $x$, then on the left-hand side of (10) there
would be other summands $\mu_{x_k}\,
K_m({x_k};a;q)\,K_{m'}({x_k};a;q)$, corresponding to these
additional points. Let us show that these additional summands do
not appear. To this end we set $m=m'=0$ in the relation (10) with
the additional summands. Since $K_0(x;a;q)=1$, we have the
equality
$$
\sum_{n=0}^\infty \frac{a^nq^{n(n+1)/2}}{(q;q)_n} + \sum_k
\mu_{x_k} =(-aq;q)_\infty  .
$$
According to the formula for the $q$-exponential function $E_q(a)$
(see formula (II.2) of Appendix II in [7]), we have
$\sum_{n=0}^\infty a^nq^{n(n+1)/2}(q;q)^{-1}_n =(-aq;q)_\infty$.
Hence, $\sum_k \mu_{x_k} =0$ and all $\mu_{x_k}$ disappear. This
means that additional summands do not appear in (10) and it does
represent the orthogonality relation for alternative $q$-Charlier
polynomials.

As we have shown, the orthogonality relation for the alternative
$q$-Charlier polynomials is given by formula (10). Due to this
orthogonality, we arrive at the following statement:
\medskip

{\bf Proposition.} {\it The spectrum of the operator $I_1$
coincides with the set of points $q^{n}$, $n=0,1,2,\cdots$. The
spectrum is simple and has one accumulation point at 0.}
\bigskip

\noindent {\bf 4. Dual alternative $q$-Charlier polynomials}
\bigskip

Now we consider the second identity in (9), which gives the
orthogonality relation for the matrix elements ${\hat a}_{mn}$,
considered as functions of $m$. Up to multiplicative factors these
functions coincide with the functions
$$
F_n(x;a|q)={}_2\phi_1 (x,-a/x;\; 0;\; q,q^{n+1}),  \eqno (11)
$$
considered on the set $x\in \{ q^{-m}\, |\, m=0,1,2,\cdots \}$.
Consequently,
$$
{\hat a}_{mn}= \left( \frac{a^nq^{n(n+1)/2}}{(q;q)_n}
\frac{(1+aq^{2m})}{(-aq^m;q)_\infty(q;q)_m\,a^m q^{m(m+1)/2}}
\right)^{1/2}   F_n(q^{-m} ;a|q)
$$
and the second identity in (9) gives the orthogonality relation
for $F_n(q^{-m} ;a|q)$:
$$
\sum_{m=0}^\infty \frac{(1+aq^{2m})}{a^m(-aq^m;q)_\infty (q;q)_m
q^{m(m+1)/2}} F_n(q^{-m};a|q)F_{n'}(q^{-m};a|q) =\frac{
(q;q)_n}{a^{n}q^{n(n+1)/2}} \delta_{nn'}. \eqno (12)
$$

The functions $F_n(x;a,b|q)$ can be represented in another form.
Indeed, taking in the relation (III.8) of Appendix III in [7] the
limit $c\to \infty$, one derives the relation
$$
{}_2\phi_1 (q^{-m},-aq^m;\; 0;\; q,q^{n+1})=(-a)^mq^{m^2}
{}_3\phi_0 (q^{-m},-aq^m, q^{-n}\; -\, ;\; q,-q^{n}/a) .
$$
Therefore, we have
$$
F_{n}(q^{-m} ;a|q)=(-a)^mq^{m^2}{}_3\phi_0
(q^{-m},-aq^{m},q^{-n};\; -\, ; \; q,-q^n/a) .
 \eqno (13)
$$
The basic hypergeometric function ${}_3\phi_0$ in (13) is a
polynomial of degree $n$ in the variable $\mu(m): = q^{-m}-
a\,q^{m}$, which represents a $q$-quadratic lattice; we denote it
by
$$
d_n(\mu (m); a;q):= {}_3\phi_0(q^{-m},-a\,q^{m},q^{-n};\; -\, ;\;
q,-q^n/a)\,.  \eqno (14)
$$
Then formula (12) yields the orthogonality relation
$$
\sum_{m=0}^\infty \frac{(1+aq^{2m})a^m}{(-aq^m;q)_\infty
(q;q)_m}\, q^{m(3m-1)/2} d_n(\mu(m)) d_{n'}(\mu(m))=\frac{
(q;q)_n}{a^nq^{n(n+1)/2}}\, \delta_{nn'}  \eqno (15)
$$
for the polynomials (14) when $a>0$. As far as we know this
orthogonality relation is new. We call the polynomials $d_n(\mu
(m); a;q)$ {\it dual alternative $q$-Charlier polynomials}. Thus,
we proved the following theorem.
 \medskip

{\bf Theorem.} {\it The polynomials $d_n(\mu (m); a;q)$, given by
formula (14), are orthogonal on the set of points $\mu(m): =
q^{-m}- a\,q^{m}$, $m=0,1,2,\cdots$, and the orthogonality
relation is given by formula (15).}
 \medskip

{\it Remark.} The duality of polynomials is the well-known notion
(see [7] and [11]) and, in particular, in the case of polynomials,
orthogonal with respect to a finite number of discrete points, it
reflects the simple fact that a finite-dimensional matrix,
orthogonal by rows, is also orthogonal by its columns ({\it cf}
(9)). There is also an analytical way of deriving a dual
orthogonality for polynomials, whose weight functions are
supported on an infinite number of discrete points (see, for
example, [12]-[14]). But this derivation is given in terms of a
dual set of functions (for instance, their explicit form (11) for
the case of alternative $q$-Charlier polynomials is given above)
and one needs to make one step further in order to extract an
appropriate family of dual polynomials from these functions.
Observe that in our approach to the duality of $q$-polynomials it
is not assumed that an initial family is orthogonal, this
orthogonality property is straightforwardly derived. Besides, one
naturally extracts from a dual set of orthogonal functions an
appropriate dual family of $q$-polynomials.
 \medskip

Let ${\frak l}^2$ be the Hilbert space of functions on the set
$m=0,1,2,\cdots$ with the scalar product
$$
\langle f_1,f_2\rangle = \sum _{m=0}^\infty\,
\frac{(1+aq^{2m})a^m}{(-aq^m;q)_\infty(q;q)_m}\, q^{m(3m-1)/2}
\,f_1(m)\,\overline{f_2(m)} ,   \eqno (16)
$$
where the weight function is taken from (15). The polynomials (14)
are in one-to-one correspondence with the columns of the
orthogonal matrix $({\hat a}_{mn})$ and the orthogonality relation
(15) is equivalent to the orthogonality of these columns. Due to
(9) the columns of the matrix $({\hat a}_{mn})$ form an
orthonormal basis in the Hilbert space of sequences ${\bf a}=\{
a_n\, |\, n=0,1,2,\cdots \}$ with the scalar product $\langle {\bf
a},{\bf a}'\rangle=\sum_n a_n{\overline{a'_n}}$. This scalar
product is equivalent to the scalar product (16) for the
polynomials $d_n(\mu (m);a;q)$. For this reason, the set of
polynomials $d_n(\mu (m);a;q)$, $n=0,1,2,\cdots$, form an
orthogonal basis in the Hilbert space ${\frak l}^2$. This means
that {\it the point measure in (15) is extremal for the dual
alternative $q$-Charlier polynomials} $d_n(\mu (m);a;q)$ (for the
definition of an extremal orthogonality measure see, for example,
in [11]). This means that the dual alternative $q$-Charlier
polynomials (14) form a complete system in the $L^2$-space with
respect to the point measure in (15). Observe that the
completeness of alternative $q$-Charlier polynomials in the
$L^2$-space with respect to the point measure in (10) is a
consequence of the fact that the operator (1) is bounded.

A recurrence relation for the polynomials $d_n(\mu (m);a;q)$ is
derived from (4). It has the form
$$
(q^{-m}- a q^{m})d_n(\mu (m))= - a d_{n+1}(\mu (m)) +
q^{-n}d_{n}(\mu (m)) - q^{-n}(1-q^{n}) d_{n-1}(\mu (m)), \eqno
(17)
$$
where $d_{n}(\mu (m))\equiv d_{n}(\mu(m);a;q)$. A $q$-difference
equation for $d_{n}(\mu(m);a;q)$ can be obtained from the
three-term recurrence relation for alternative $q$-Charlier
polynomials.

Note that for the polynomials $d_n(\mu (m);a;q^{-1})$ with $q<1$
we have the expression
$$
d_n(\mu (m);a;q^{-1})={}_3\phi_2(q^{-m},-a\,q^{m},q^{-n};\; 0,0\,
;\; q,q) .  \eqno (18)
$$
However, the recurrence relation for these polynomials (which can
be obtained from the relation (17), does not satisfy the
positivity condition $A_nC_{n+1}>0$, that is, they are not
orthogonal polynomials for $a>0$ (as it is the case for
alternative $q$-Charlier polynomials). This positivity condition
holds only if we demand $a<0$. In this case, the polynomials (18)
are the continuous big $q$-Hermite polynomials $H_n(x;a|q)$ (for
an explicit form of these polynomials see, for example, [1],
formula (3.18.1)), which are orthogonal on a certain continuous
set of points.
\bigskip

\noindent{\bf Acknowledgments}

\medskip

This research has been supported in part by the SEP-CONACYT
project 41051-F and the DGAPA-UNAM project IN112300 "\'Optica
Matem\'atica". A.~U.~Klimyk acknowledges the Consejo Nacional de
Ciencia y Technolog\'{\i}a (M\'exico) for a C\'atedra Patrimonial
Nivel II.

\end{document}